\documentclass[a4paper,11pt]{amsart}

\usepackage{url}
\usepackage{color}
\definecolor{darkgreen}{rgb}{0,0.5,0}
\usepackage[
        colorlinks, citecolor=darkgreen,
        backref,
        pdfauthor={V. Arul, J. Steffen Mueller},
        pdftitle={Rational points on X_0^+(125)$}
]{hyperref}

\usepackage{amssymb,amsmath, amsthm}
\usepackage{comment}
\usepackage{algorithmic,algorithm}
\usepackage[OT2,OT1]{fontenc}
\usepackage{graphicx}
\usepackage{epsfig}
\usepackage{fancyhdr}
\usepackage{enumerate}
\usepackage{colonequals}

\usepackage[alphabetic,backrefs]{amsrefs} 
\usepackage[all]{xy}

\newcommand\cyr{%
  \renewcommand\rmdefault{cmr}%
  \renewcommand\sfdefault{wncyss}%
  \renewcommand\encodingdefault{OT2}%
  \normalfont\selectfont}
\DeclareTextFontCommand{\textcyr}{\cyr}

\definecolor{red}{rgb}{0.9,0,0}
\definecolor{purple}{rgb}{0.8,0,0.6}

\numberwithin{equation}{section}

\newtheorem{thm}{Theorem}[section]

\newtheorem{prop}[thm]{Proposition}

\newtheorem{conj}[thm]{Conjecture}

\theoremstyle{definition}

\theoremstyle{remark}

\newcommand\Q{\mathbb{Q}}

\newcommand\F{\mathbb{F}}

\newcommand\Z{\mathbb{Z}}

\newcommand\NS{\mathop{\rm NS}\nolimits}

\newcommand{\et}{\scriptstyle \mathrm{\acute{e}t}}

\newcommand{\rank}{\operatorname{rk}}

\newcommand{\GL}{\operatorname{GL}}

\AtBeginDocument{%
     \def\MR#1{}
}

\begin{document}

\title{Rational points on $X^+_0(125)$}
\author{Vishal Arul}
\email{varul.math@gmail.com}
\address{ Vishal Arul,
  Department of Mathematics, University College London, United Kingdom}

\author{J. Steffen M\"uller}
\email{steffen.muller@rug.nl}
\address{ J. Steffen M\"uller,
  Bernoulli Institute, 
  University of Groningen,
  Nijenborgh 9,
  9747 AG Groningen,
  The Netherlands
}

\date{\today}


\begin{abstract} \setlength{\parskip}{1ex} \setlength{\parindent}{0mm}
  We compute the rational points on the Atkin--Lehner quotient
  $X^+_0(125)$ using the quadratic Chabauty method. Our work completes the study of
  exceptional rational points on the curves $X^+_0(N)$ of
  genus between 2 and 6.
  Together with the work of several authors, this completes the proof of a conjecture of
  Galbraith.
\end{abstract}

\maketitle




\section{Introduction}\label{S:intro}

  For an integer $N>1$, let $w_N$ denote the Atkin--Lehner involution on the modular curve $X_0(N)$. Then the
  {\em Atkin--Lehner quotient of level $N$} 
  \begin{equation*}\label{}
    X^+_0(N)\colonequals X_0(N)/\langle w_N\rangle 
  \end{equation*}
  is a smooth projective curve over $\Q$.
  A rational point on $X^+_0(N)$
  lifts to a  quadratic point on $X_0(N)$, which provides motivation for the computation of
  $X_0^+(N)(\Q)$; see for
  instance~\cite{Box21}.

  A non-cuspidal point on $X^+_0(N)$ corresponds
  to an unordered pair $\{E,E'\}$ of elliptic curves, together with a cyclic isogeny $E\to E'$ of
  degree $N$. 
  In particular, a non-cuspidal point in $X^+_0(N)(\Q)$ such that the associated 
  elliptic curves do not have complex multiplication corresponds to a quadratic
  $\Q$-curve~\cite{Elk04} without complex multiplication. 
Following
  Galbraith~\cite{Gal02}, we call such a rational point
  \textit{exceptional}.
  It turns out that exceptional rational points are quite rare. 
    Guided by extensive numerical computations, Galbraith made the following conjecture
  in~\cite{Gal02}.
  We write $g^+_N$  for the genus of $X^+_0(N)$.

  \begin{conj}\label{ConjGal}
    Let $2\le g^+_N\le 5$. Then $X^+_0(N)$ contains exceptional rational points if and only
    if 
    \begin{equation*}\label{}
  N \in \{73, 91,103,125,137, 191,311\}\,.
    \end{equation*}
  \end{conj}
The main result of this short note is the following modest contribution to the study of the rational
points on $X^+_0(N)$.
\begin{prop}\label{P:125}
  There are precisely~6 rational points on $X_0^+(125)$. One of these is a cusp, four are
  CM points 
with respective discriminants $D=-19,-16,-11,-4$, and one is
exceptional, corresponding to a quadratic $\Q$-curve with $j$-invariant
satisfying
  \begin{align*}
    11^5 j=  &-2140988208276499951039156514868631437312\\&\pm
  94897633897841092841200334676012564480\sqrt{509}\,.
  \end{align*}
\end{prop}
These points, the discriminants and the $j$-invariant of the exceptional
rational point  were already found by Galbraith~\cite[Section~10]{Gal02}. To show that there
are no other rational points on this curve, we use the quadratic Chabauty method~\cite{BD18, BD20, BDMTV19, EL21,
BDMTV2, BMS}, combined with the Mordell--Weil sieve~\cite{BS10}, as explained in~\cite{BBBLMTV19,
BDMTV2}.
The quadratic Chabauty method makes Kim's non-abelian Chabauty program~\cite{Kim05, Kim09} explicit in the
simplest non-abelian case using $p$-adic heights.

It was already known by work of Momose~\cite{Mom87},
Galbraith~\cite{Gal02} and Arai--Momose \cite{AM10} that there are no composite integers
$N$ such that $2\le g_N^+\le 6$ and $X^+_0(N)$ has an exceptional rational point, except
for $N=91,125$ (where an exceptional point was found by Galbraith) and
possibly for $N=169$.
See Section~\ref{S:comp} below. The
rational points (in fact the $\Q(i)$-points) on $X_0^+(91)$  were computed by 
Balakrishnan, Besser, Bianchi, and the second-named author, see~\cite[Example~7.1]{BBBM21}.
Balakrishnan et al.\ showed in~\cite{BDMTV19} that $X_0^+(169)$
contains no exceptional rational points. For both of these computations, the quadratic
Chabauty method was used. 
We obtain:
\begin{thm}\label{T:main}
  The only composite levels $N$ such that $X^+_0(N)$ has genus $g^+_N\in\{2,\ldots,6\}$
    and contains an exceptional rational point are $N=91$ and $N=125$. 
    There are two exceptional rational points on $X^+_0(91)$, with
    $j$-invariants satisfying
    \begin{equation*}
      2^{14}j =
      -270483906936119115236875\pm6098504215856136863625\sqrt{-87}
    \end{equation*}
    and 
      \begin{equation*}
        2^{92}j =
      -8366877442964720618049886816125\pm32028251460268098916979319375\sqrt{-87}.
      \end{equation*}
    The $j$-invariant of the exceptional rational 
    point on $X^+_0(125)$ is given in
    Proposition~\ref{P:125}.
\end{thm}
The $j$-invariants were already computed by Galbraith~\cite{Gal02}.

For prime levels $N$ such that $2\le g_N^+\le 6$, the rational points on $X^+_0(N)$ 
were computed in~\cite{BBBLMTV19, BDMTV2,AABCCKW}. More precisely:
\begin{itemize}
  \item Balakrishnan et al.~\cite{BBBLMTV19} computed $X^+_0(N)(\Q)$ for $N \in
    \{67,73,103\}$.
  \item Balakrishnan et al.~\cite{BDMTV2} computed $X^+_0(N)(\Q)$ for $N \in
    \{107, 167, 191\}$ (the remaining genus~2 prime levels) and for all prime $N$ such that
    $g_N^+=3$.
  \item Adžaga et al.~\cite{AABCCKW} computed $X^+_0(N)(\Q)$ for all prime
    $N$ such that $g_N^+\in\{4,5,6\}$.
\end{itemize}
All exceptional points that occurred in these computations had already been found by
Galbraith.
Together with Theorem~\ref{T:main}, this implies:
\begin{thm}\label{T:}
  Conjecture~\ref{ConjGal} holds.
\end{thm}
Of course, the upper bound~5 on the genus in Galbraith's conjecture seems rather arbitrary
and is likely due to computational limitations. 
According to Elkies~\cite[p. 44]{Elk98}, one would expect that there are no
exceptional rational points on $X_0^+(N)$ for $N\gg 0$. The explicit methods of this paper
are not suitable to prove such a statement. However, there are some known results that hold
for infinitely many levels.
For instance, Dogra and Le Fourn~\cite{DLF21} have
shown that the Chabauty--Kim set $X_0^+(N)(\Q_p)_2$ is finite for all prime $N$ such that
$g^+_N>1$. Furthermore, the results of Momose~\cite{Mom86, Mom87} and
Arai-Momose~\cite{AM10} apply to infinitely many composite $N$.

It would also be interesting to compute the rational points on quotients of
$X_0(N)$ by more general groups of Atkin--Lehner involutions. 
  Adžaga, Chidambaram, Keller and Padurariu~\cite{ACKP22} have recently completed the computation of the
  rational points on the quotient $X_0^*(N)$ by the full Atkin--Lehner subgroup whenever
  this curve is hyperelliptic of genus $>1$. See also~\cite{BGX20}.
\subsection*{Acknowledgements}  
We thank Timo Keller for comments on an earlier version of this paper, Jennifer Balakrishnan, Netan Dogra,
Nikola Adžaga, Lea Beneish, Mingjie Chen, Shiva Chidambaram,
Timo Keller and Boya Wen for helpful discussions and the anonymous referee
for useful suggestions.
SM was supported by DFG grant MU 4110/1-1 and by an NWO Vidi grant. 

\section{Composite level}\label{S:comp}
  We first list all composite levels $N$ such that 
  $2\le g^+_N\le 6$:
  \begin{align*}
    g^+_N=2: \quad N = &\,42, 46, 52, 57, 62, 68, 69, 72, 74, 77, 80, 87, 91, 98, 111, 121, 125, 143\\
    g^+_N=3:\quad N = &\,58, 60, 66, 76, 85, 86, 96, 99, 100, 104, 128, 169\\
    g^+_N=4:\quad N = &\,70, 82, 84, 88, 90, 92, 93, 94, 108, 115, 116, 117, 129, 135,
    147, 155, 159,\\ 
        &\,161, 215\\
    g^+_N=5:\quad N = &\, 78, 105, 106, 110, 112, 122, 123, 133, 134, 144, 145, 146, 171,
    175, 185, 209\\
    g^+_N=6:\quad N = &\, 118, 124, 136, 141, 152, 153, 163, 164, 183, 197, 203, 211,
    221,223, 269, \\&\, 271, 299, 359
  \end{align*}
The list is taken from Table~2 of~\cite{AABCCKW}, and follows from 
a lower bound for $g_N^+$
 in terms of $N$ given in Proposition~4.4 of~\cite{AABCCKW}. This bound was obtained by
 expressing $g_N^+$ in terms of the genus of
$X_0(N)$ and the class number of certain quadratic orders. Bounds on the former are
classical, and the latter can be bounded using the Dirichlet class number formula.

  We now discuss previous work on the computation of $X_0^+(N)(\Q)$ for composite $N$.
  In~\cite{Mom87}, Momose showed that there are no rational points on $X_0^+(N)$ if
  $N$ is composite and contains a sufficiently large prime factor, in the following sense.
  \begin{thm}\label{T:Mom87}(Momose,~\cite[Theorem~0.1, Proposition~2.11]{Mom87})
    The curve $X_0^+(N)$ does not contain an
  exceptional rational point when $N$ is composite and one of the following conditions
  holds:
\begin{enumerate}
    \item There is a prime $p\mid N$ such that $p\ge 11$, $p\notin\{13,37\}$ and the
      isogeny factor $J_0^-(p)$ of $J_0(p)$ has Mordell--Weil rank~0 over $\Q$. The latter holds
      for $p=11$ and all primes $p\in \{17,\ldots,300\}-\{151,199,227,277\}$.
    \item $g^+_N\ge 1$ and at least one of the following: $26\mid N$,
      $27\mid N$, $35\mid N$.
    \item $g^+_N\ge 1$ and $49\mid N$; moreover
      $N/49$ 
      \begin{itemize}
        \item is divisible by 7 or 9, or 
        \item is divisible by a prime $q\equiv 2 \bmod{3}$ or
        \item is not divisible by 7 but satisfies $\left(\frac{-7}{m}\right)=-1$.
      \end{itemize}
  \end{enumerate}
  \end{thm}
Momose's techniques do not apply to multiples of $37$. This case was resolved much later by Arai and Momose.
  \begin{thm}\label{T:AM10}(Arai-Momose,~\cite[Theorem~1.2]{AM10})\
    The curve $X_0^+(N)$ does not contain an
  exceptional rational point when $N$ is composite and divisible by~37.
  \end{thm}
  
The only composite levels $N$ such that $2\le g^+_N\le 6$ and such that $N$ is not covered
by Theorem~\ref{T:Mom87} or Theorem~\ref{T:AM10} are
 \begin{align*}
    g^+_N=2:\quad& 42, 72, 80, 91, 125 \\
g^+_N=3:\quad& 60, 96, 100, 128, 169\\
g^+_N=4:\quad& 84, 90, 117\\ 
g^+_N=5:\quad&  112, 144\\
   g^+_N=6:\quad&  \textrm{none}\,.
\end{align*}
For all of these, Galbraith constructs explicit models in~\cite{Gal02} using the techniques
of~\cite{Gal99}. For $N\notin\{91,125,169\}$, he also computes the
rational points using a morphism to an elliptic curve with rank~0. Galbraith gives no details,
but it is not hard to find such a morphism, for instance by computing the
automorphisms of the curve. We provide code to check these computations
at~\cite{AMCode}.
Note that for some of these levels, such as $N=128$,~\cite[Theorem~0.1]{Mom86} also implies that
there are no exceptional rational points.

\section{The remaining cases}\label{S:QC}
In~\cite{Gal02}, Galbraith computes the rational points of small height on the curves
$X_0^+(N)$ for $N\in \{91,125,169\}$
 and conjectures that these are all rational points.
The techniques mentioned at the end of Section~\ref{S:comp} cannot be used
here, because these curves do not cover a curve with finitely many rational
points. Moreover, they have Mordell--Weil rank equal to the genus, so that the method of Chabauty
and Coleman~\cite{McP12} is not
applicable.  

Fortunately, all three curves have Jacobians with real multiplication, and one may apply the 
quadratic Chabauty method~\cite{BD18, BDMTV19, BD20} to find a finite set of $p$-adic
points containing the rational points. When the genus is~2, one typically needs to combine
this 
with the Mordell--Weil sieve~\cite{BS10} as
in~\cite{BBBLMTV19,BDMTV2},
to identify the rational points among the solutions. We
describe this computation for $N=125$ in Section~\ref{S:125} below.

The curve $X_0^+(91)$ is a bielliptic genus~2 curve of Mordell--Weil rank~2. Its label in
the {\tt
LMFDB}~\cite{lmfdb} is
{\href{https://www.lmfdb.org/Genus2Curve/Q/8281/a/8281/1}{8281.a.8281.1}
and it can be described by the equation $y^2 = x^6-3x^4+19x^2-1$.
Hence, the explicit methods developed by Balakrishnan and Dogra in~\cite[Section~8]{BD18}
apply. 
In fact, the $\Q(i)$-points on this curve were found in~\cite[Example~7.1]{BBBM21}
using an extension of these techniques to curves over number fields and the
Mordell--Weil
sieve.
This shows that that $\#X_0^+(91)(\Q)=10$, as
predicted by Galbraith. 

The curve $X^+_0(169)$ is a non-hyperelliptic curve of genus~3. It is isomorphic to both
the split Cartan modular curve $X_{\mathrm{s}}^+(13)$ and the nonsplit Cartan modular
curve $X_{\mathrm{ns}}^+(13)$. In~\cite{BDMTV19} the quadratic Chabauty method was
applied to show that $X^+_0(169)$ has precisely~7 rational points, as predicted by
Galbraith. There was no need for the Mordell--Weil sieve because it was
possible to construct two $17$-adic functions whose common zero set is precisely the set
of rational points (it turns out that one can actually work $3$-adically,
see~\cite{QCCode}).

\section{Quadratic Chabauty for $X^+_0(125)$}\label{S:125}
The final curve to consider is $X\colonequals X^+_0(125)$, with {\tt LMFDB} label \href{https://www.lmfdb.org/Genus2Curve/Q/15625/a/15625/1}
{15625.a.15625.1}.
Galbraith~\cite[Section~10]{Gal02} finds the model
\begin{equation}\label{125model}
  y^2 = x^6+2x^5+5x^4+10x^3+10x^2+8x+1\,,
\end{equation}
and the six rational points $\infty_{\pm}, (0,\pm1), (-2,\pm 5)$. Of these,
$\infty_-, (0,1), (0,-1), (-2, 5)$ are Heegner points,
$\infty_+$ is a cusp, and $(-2,-5)$ is
exceptional.

Since the Galois group of the polynomial $x^6+2x^5+5x^4+10x^3+10x^2+8x+1$ is non-abelian
of order~60, it is not feasible to apply elliptic curve Chabauty, as introduced by Bruin~\cite{Bru03} and used, for
instance, to find the rational points on many genus~2 Atkin--Lehner quotients $X_0^*(N)$
by Bars, Gonz\'alez and Xarles~\cite{BGX20}.

We compute the rational points on $X$ by applying quadratic Chabauty for $p=29$. Our code~\cite{AMCode}
is written in {\tt
Magma}~\cite{BCP97}, and is based on the package {\tt QCMod}, available from~\cite{QCCode}. We describe the computation below,
referring to~\cite{BD18, BD20} for the theoretical background and to \cite{BDMTV19,
BBBLMTV19, BDMTV2} for details on the explicit methods. Alternative approaches to the
quadratic Chabauty method have been introduced by Edixhoven and Lido~\cite{EL21} and by Besser,
Srinivasan and the second-named author~\cite{BMS}.

The quadratic Chabauty method for rational points, described explicitly in~\cite{BDMTV2},
requires (in its simplest form) that
$\rank J(\Q)=g$, that $\rank \NS(J)>1$ and that the closure of $J(\Q)$ in $J(\Q_p)$ has
finite index, where $p$ is a prime of good reduction for $X$. As discussed
above, the first two conditions are satisfied for $X$. The third
condition means that the method of Chabauty--Coleman is not applicable; it is satisfied for
all Atkin--Lehner quotients $X^+_0(N)$ and all $p$ by~\cite[Lemma~7]{DLF21}.

The quadratic Chabauty method also requires an explicit description of the semi-stable
reduction. Clearly our curve $X$ has good reduction away from~5.
For many modular curves of level $N$, the semi-stable reduction at primes $\ell\mid N$ is
known.
For instance, 
when $\ell$ is a prime and $N=\ell$ or $N=\ell^2$, Edixhoven~\cite{Edi89, Edi91} has described the semi-stable model of $X_0(N)$ and the action of
$w_N$ on it, from which one may deduce the semi-stable reduction of $X^+_0(N)$. More
recently, Edixhoven and Parent~\cite{EP21} found the semi-stable reduction of the modular curve
associated to any maximal subgroup of $\GL_2(\F_\ell)$.
However, we are not aware of any statements in the literature that would allow us to find the
semi-stable reduction of $X^+_0(125)$ at $5$ explicitly using a modular approach, though the
techniques of Weinstein~\cite{Wei16} may be of use here.
Instead, we apply Qing Liu's code
\href{https://www.math.u-bordeaux.fr/~qliu/G2R/index.html}{genus2reduction}
(for instance through
\href{http://pari.math.u-bordeaux1.fr/dochtml/html-stable/Elliptic_curves.html#genus2red}{Pari/GP} or
\href{https://doc.sagemath.org/html/en/reference/arithmetic_curves/sage/interfaces/genus2reduction.html}{SageMath}).
This shows that $X^+_0(125)$ has potentially good reduction
at~5.  
While potentially good reduction is not strictly necessary (see~\cite[\S{5.4}]{BDMTV2}) for
the quadratic Chabauty method, it greatly simplifies its 
application in practice; see~\cite[\S{3.1}]{BDMTV2}
and~\cite[Section~12]{BeD19} for a discussion.

We use the prime $p=29$ of good reduction and the model
\begin{equation}\label{125model2}
  y^2= -1487x^6 - 3238x^5 - 2915x^4 - 1390x^3 
      - 370x^2 - 52x - 3
\end{equation}
of $X$.
While this has larger coefficients than~\eqref{125model}, it has the advantage that there
are no $\F_{29}$-rational points at infinity; hence we have $X(\Q) = Y(\Q)$,
where $Y$ is the affine
plane curve defined by~\eqref{125model2}.
 Because there are also no Weierstrass points 
in $X(\F_{29})$, all rational points are contained in non-Weierstrass
residue disks of $Y(\Q_{29})$.
We fix the base-point 
$b =\left(-\frac{1}{2} ,\frac{ 1 }{8}\right) \in X(\Q)$ for the Abel-Jacobi
map $X\to J$.

We compute the Frobenius matrix $F_{29}$ on
$H^1_{\mathrm{dR}}(X_{\Q_{29}})$ using the algorithm of Tuitman~\cite{Tui16,Tui17}.
From $F_{29}$, we find the Hecke operator $T_{29}$ using Eichler-Shimura, and form a nice
correspondence $Z$ on $X\times X$ (corresponding to a nontrivial element of $\ker(\NS(J)
\to \NS(X))$) as in~\cite[\S{6.4}]{BDMTV19}. The correspondence $Z$ induces a fundamental group quotient $U=U_Z$
whose Chabauty--Kim set $X(\Q_{29})_U$ is finite
by~\cite[Lemma~3.1]{BD18}.
From this, we may find a global (respectively local)  ${29}$-adic Galois representation $A_Z(x)
= A_Z(b,x)$ with graded pieces $\Q_{29}, H^1_{\et}(\bar{X},
\Q_p)^\vee$ and
$\Q_{29}(1)$  for each $x\in
X(\Q)$ (respectively $x\in X(\Q_{29})$), see~\cite[Section~5]{BD18} and~\cite[\S3.4]{BDMTV19}.
By~\cite[Section~5]{BD18}, the local height $x\mapsto h_{29}(A_Z(x))$ constructed by
Nekov{\'a}{\v{r}}~\cite{Nek93} is a locally analytic
function on all of $X(\Q_{29})$, and 
the global ${29}$-adic height $h(A_Z(x))$ extends to a locally analytic function $h\colon
X(\Q_{29})\to \Q_{29}$ as well.
Since our curve has potentially good reduction, we have $h(x)
=h_{29}(A_Z(x))$ for $x\in X(\Q)$ by~\cite[Lemma~3.2]{BDMTV19}. This
equality only holds for finitely many points in $X(\Q_{29})$, so we obtain  
a locally analytic {\em quadratic Chabauty function}
\begin{equation*}\label{}
  \rho\colon X(\Q_{29})\to \Q_{29}\,;\quad x\mapsto h(x)-h_{29}(A_Z(x))
\end{equation*}
that vanishes along $X(\Q)$ and has only finitely many zeros. 

On the model~\eqref{125model2}, the~6 small rational points found by Galbraith are:
\begin{equation*}\label{}
  X(\Q)_{\mathrm{known}} = \left\{ \left(-\frac{1}{2} ,\pm\frac{ 1 }{8}\right), 
  \left(-\frac{1}{3} ,\pm\frac{1}{27} \right), 
  \left(-\frac{1}{4},\pm\frac{5}{64}\right)\right\}\subset Y(\Q)\,.
\end{equation*}
The curve $X$ has sufficiently many rational points in the sense of~\cite[\S{3.3}]{BDMTV2} so that we may determine the global
height pairing $h$ as a bilinear pairing on the tangent space
$\mathrm{H}^0(X_{\Q_{29}},\Omega^1)^{\vee}$ using the values $h(A_Z(x))$ for $x\in
X(\Q)_{\mathrm{known}}$. In our example, one needs three rational points with pairwise distinct $x$-coordinates. This suffices to extend $x\mapsto h(A_Z(x))$ to a locally
analytic function $h\colon X(\Q_{29})\to \Q_{29}$, see~\cite[\S{3.3}]{BDMTV2}.

For a point $x\in X(\Q_{29})$, the local height~$h_{29}(A_Z(x))$ 
can be expressed in  terms of the
filtered $\phi$-module $\mathrm{D}_{\mathrm{cris}}(A_Z(x))$, where
$\mathrm{D}_{\mathrm{cris}}$ is Fontaine's functor.
While $x\mapsto h_{29}(A_Z(x))$ is locally analytic on all of $X(\Q_{29})$, the
explicit methods for its computation discussed in~\cite{BDMTV19, BDMTV2} require
that the Frobenius lift of Tuitman~\cite{Tui16, Tui17} is defined at $x$.
This holds if $x$ lies in an affine non-Weierstrass disk; 
by construction, all points in $Y(\Q_{29})$ and in $X(\Q)$
are contained in such a disk.
Using the techniques of~\cite[Sections 4, 5]{BDMTV19}, we may compute the Hodge filtration and the Frobenius
structure of $\mathrm{D}_{\mathrm{cris}}(A_Z(x))$ for  $x \in Y(\Q_{29})$, and we obtain an
expansion of $\rho$ as a power series on every residue disk in $Y(\Q_{29})$.

We compute that the function $\rho$ indeed vanishes along
$X(\Q)_{\mathrm{known}}$; it  also vanishes in 22 additional points $x\in Y(\Q_{29})-
X(\Q)_{\mathrm{known}}$. The points are provably correct to $O(29^5)$,
since the code incorporates the precision analysis in~\cite[\S4]{BDMTV2}.
Following~\cite[\S{3.4}]{BDMTV19}, we apply the Mordell--Weil sieve with the prime $v= 1399$ to show that none of these are rational.
More precisely, suppose that $x\in Y(\Q_{29})- X(\Q)_{\mathrm{known}}$ were rational. Then
we could write $[x-b] = a_1P_1+a_2P_2$ with $a_1,a_2\in \Z$, where $J(\Q)=\langle P_1,P_2\rangle\cong\Z^2$. Under
this assumption, we find $a_i\bmod{29}$ from $x$ for $i\in \{1,2\}$ using linearity of
Coleman integrals of holomorphic differentials and we 
compute the image of the putative point $[x-b]$ in $J(\F_{1399})\cong (\Z/(29\cdot
50))^2$. We find that none of the images of the points $x \in Y(\Q_{29})-
X(\Q)_{\mathrm{known}}$ are in the 
image of
$X(\F_{1399})$ under our Abel--Jacobi map, so we derive a contradiction.
Since $X(\Q)\subset Y(\Q)$,
This proves that we indeed  have
\begin{equation*}\label{}
  X(\Q) = X(\Q)_{\mathrm{known}}\,,
\end{equation*}
which completes the proof of Theorem~\ref{T:main}.
\bibliographystyle{alpha}
\bibliography{References}
\end{document}